# Interlocking of convex polyhedra: towards a geometric theory of fragmented solids

A. J. Kanel-Belov[§], A. V. Dyskin[*], Y. Estrin[†], E. Pasternak[**], I. A. Ivanov-Pogodaev[+]

[§] Moscow Institute of Open Education, B.Vlasievskii-11, Moscow, Russia, International University Bremen, Campus ring 1, 28759, Bremen, Germany. E-mail: *kanel@mccme.ru*

[*] School of Civil and Resource Engineering, The University of Western Australia, 35 Stirling Highway, Crawley WA 6009, Australia, E-mail: *adyskin@cyllene.uwa.edu.au*

[†] ARC Centre of Excellence for Design in Light Metals,
Department of Materials Engineering, Monash University, Clayton, Vic. 3800, Australia and CSIRO Division of Manufacturing and Materials Technology, Clayton, Vic. 3168, Australia E-mail: *estrin@mech.uwa.edu.au*

[**]School of Mechanical Engineering, The University of Western Australia, 35 Stirling Highway, Crawley WA 6009, Australia,

[+] Department of Mechanics and Mathematics, Moscow State University, Vorobievy Gory, Moscow, 119898, Russia E-mail: *ivanov-pogodaev@mail.ru*

**Abstract**
We present structures comprised of identical convex polyhedra which are interlocked geometrically. These sets cannot be disassembled by removing individual polyhedra by translations and/or rotations. The shapes that permit interlocking arrangements include all five platonic solids. Criteria for interlocking based on transformations of the cross-sections of the elements in a 3D reconstruction of a layer from its middle cross-section are formulated. A generalization to higher dimensions is also given. In particular, an interlocking layer of four-dimensional cubes is described.

## 1 Introduction

This article presents a collection of structures of identical convex polyhedra interlocked geometrically. Methods for generating such structures are discussed along with their general geometrical properties.

Consider a set of contacting convex figures in $\mathbf{R}^2$. It is not difficult to prove that one of these figures can be moved out of the set by translation without disturbing others. Therefore, any set of planar figures can be disassembled by moving all figures *one by one* (as opposed to simple homothetic transformation of the plane, which automatically puts the figures out of contact with each other).

One may think that this statement holds for higher dimensions as well. However, attempts to generalize it to $\mathbf{R}^3$ have been unsuccessful. As a counterexample, structures consisting of convex solid bodies which cannot be disassembled by removing individual bodies were suggested in Ref. [1]. Their example was complicated and was constructed in a non regular way.

The simplest possible interlocking structures would be layer-like ones, i.e. sets of convex bodies situated in a layer between two parallel planes [2].

First, we consider two types of regular tetrahedra. Let us call them A- tetrahedron and B-tetrahedron (figure 1). The length of the upper edges of both tetrahedra, which are mutually perpendicular, is taken to be unity. Both tetrahedra are completely defined by their upper edges. It is clear that it is possible to place the tetrahedra in contact with each other in such a way that the front end of the upper edge of B and the middle of the upper edge of A become congruent. In this case the distance between the centers of the tetrahedra equals 1.

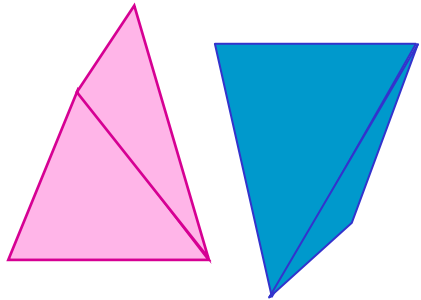

Figure 1. Tetrahedra of types A and B

Now consider the set of points on the plane such that their both coordinates are integers. Suppose that the points with both odd or both even coordinates are «white» and the points with one odd and one even coordinate are «black». Suppose that the white points are the centers of A-tetrahedra and black points are the centers of B-tetrahedra. Figure 2 shows a fragment of such a set of tetrahedra.

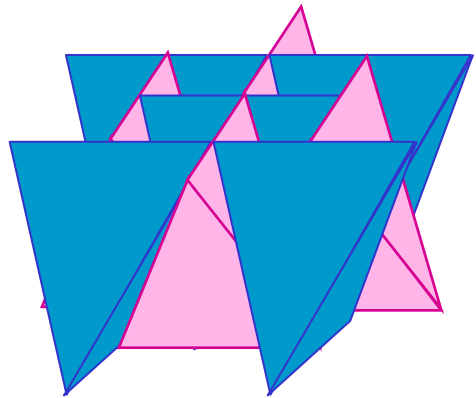

Figure 2. Fragment of interlocking set of tetrahedra.

Consider one of the A- tetrahedra, call it X. There are four neighboring B-tetrahedra: left, right, top and bottom. It is clear that four faces of X adjoin four neighboring B-tetrahedra. Furthermore, left and right neighbors prevent X from moving up while top and bottom ones prevent X from moving down. It is obvious that B-tetrahedra are also «locked» by their neighbors. This situation is called «*interlocking*».

Definition 1. Let us consider a set of solids (that may be infinite). Assume that for any solid X the following condition holds: if all the solids except X have a degree of freedom equal to 0 then X also has a degree of freedom equal to 0. In this case the set is called an *interlocking structure*. The state of *interlocking* means that no one solid can be removed from the system, since all other solids are immobile.

Below we consider various types of interlocking structures and methods of creating them.

First we consider how to create an interlocking structure using regular lattices.

Suppose that we have a regular square lattice in two colors like an infinite chessboard. Let us place an arrow at each edge using the following rules:

   1. All arrows have the same length;
   2. All arrows start in the middle of the corresponding edge and are perpendicular to it;
   3. For black squares top and bottom arrows are outside while left and right ones are inside;
   4. For white squares top and bottom arrows are inside while left and right ones are outside.

Now let us consider a plane based on an edge and lopsided to the direction of a corresponding arrow (Figure 3)

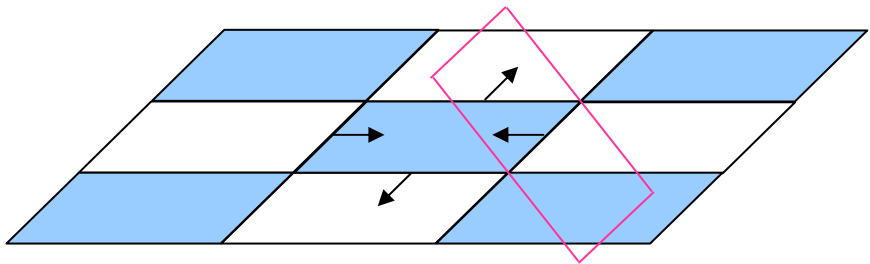

Figure 3. Lopsided plane.

Let us create such a plane for each edge using the same tilt angle α. Consider four planes of a particular square S. Obviously, these planes, if inclined at a proper angle, form a tetrahedron. This tetrahedron intersects with our initial chessboard plane by the square S.

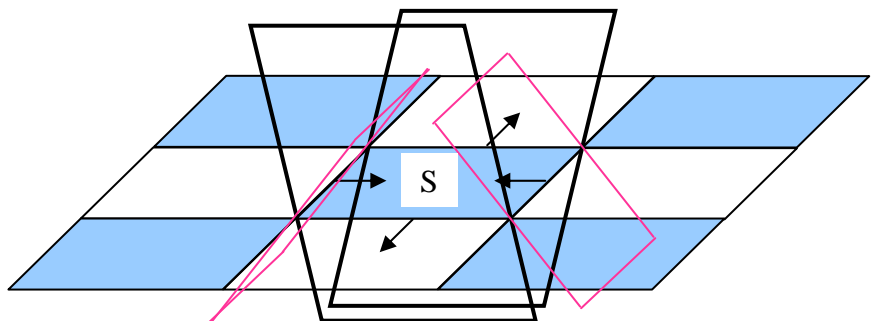

Figure 4. Lopsided planes of the square

It is clear that we obtain A-tetrahedra for the black squares and B ones for the white. Every tetrahedron is locked within the layer: The left and right neighbors prevent A- tetrahedra from moving up, while top and bottom neighbors prevent it from moving down.

Our initial chessboard plane with arrows thus becomes a diagram representing an interlocking scheme. The arrows show how the tetrahedra lock their neighbors.

This process is referred to as the «*moving cross-section procedure*».

The orientation of an arrow indicates the direction of movement of the respective face of the polygon in the cross-section upon *upward* displacement of the section plane. The length of the arrow is proportional to tan$\alpha$, the slope of the face. For the square lattice considered above, this arrow pattern is shown in Fig. 3.

The above structure provides an example of an assembly in which no element can be removed without disturbing its neighbors. It is appealing to use such arrangements in engineering as a new design principle. This principle allows one to build structures of convex elements without a binder phase or connectors, i.e., structures which (a) possess high resistance to fracture propagation (fractures getting stopped at interfaces between the elements) and free of stress concentrations associated with mechanical connectors used in conventional assemblies of building blocks [3-7]. Thus, there is a practical need to find and characterize interlocking sets. It is interesting to note that such a configuration of tetrahedra was first found by a civil engineer, Glickman [8], who proposed a pavement system based on truncated tetrahedra. This structure went unnoticed by the mathematics community, however.

The structures described are obtained by tiling the plane with equal squares. A.J. Kanel-Belov considered another regular tiling - the hexagonal (honeycomb) one and found a new type of interlocking solids, which turn out to be simple cubes! The recipe for constructing new interlocking structures based on this tiling is as follows. Put arrows of equal length assigned to the faces of a hexagon in such a way that for any hexagon within the tiling the arrows heading inwards and outwards alternate, Fig. 5. Consider a layer of hexagonal prisms with their honeycomb type middle section. Incline the lateral faces of the prisms as prescribed by the arrows in accordance with the aforementioned convention. If these inclined faces are extended as in the above example, an interlocking arrangement of cubes (Fig. 6), octahedra (Fig. 7) or dodecahedra (Fig. 8) is generated, depending upon the angle $\alpha$ [9]. All three resulting shapes correspond to the *platonic bodies*, as does the tetrahedron one considered above. It is interesting to note that for cubes the angle $\alpha = \sin^{-1}(\sqrt{3}/3)$ is the same as for tetrahedra. For octahedra this angle is smaller: $\alpha = \sin^{-1}(1/3)$.

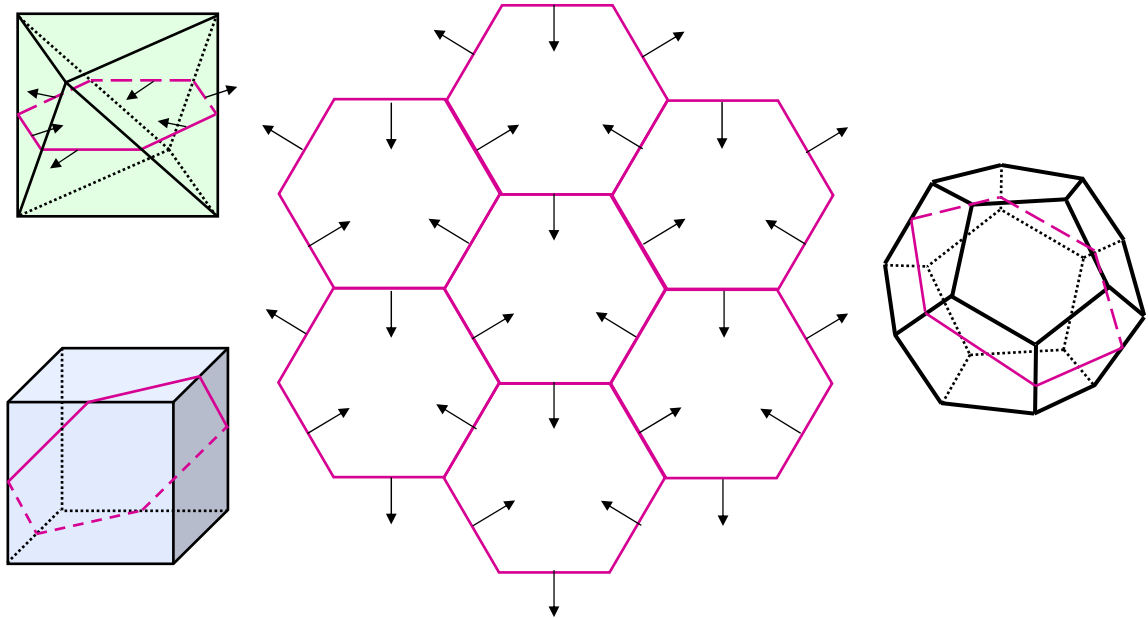

Figure 5. Hexagonal tiling of the plane and the associated platonic solids with their respective hexagonal middle sections. Orientations of arrows indicate the inclinations of faces of the modified hexagonal prisms that provide their interlocking.

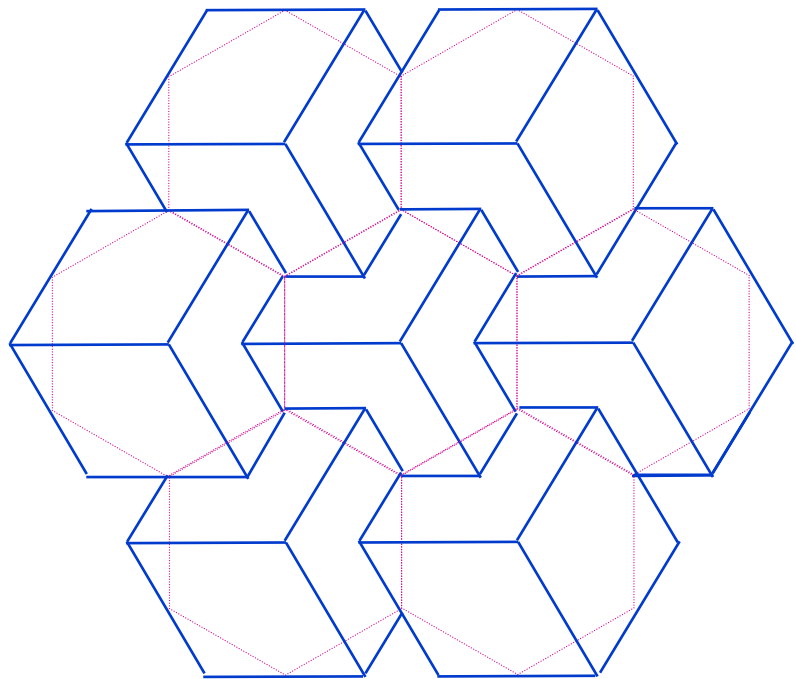

Figure 6. Fragment of interlocking set of cubes.

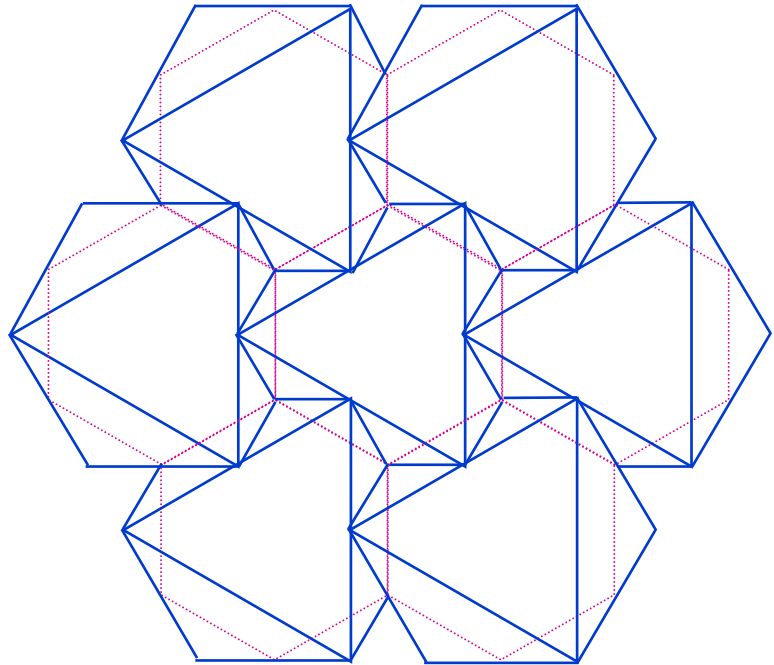

Figure 7. Fragment of interlocking set of octahedra.

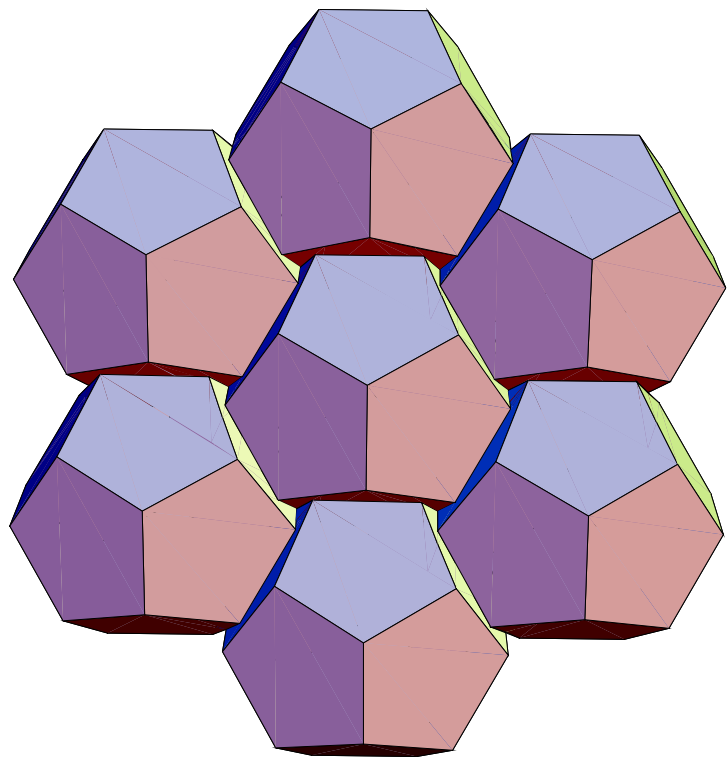

Figure 8. Fragment of interlocking set of dodecahedra normal to the symmetry axis of 3$^{rd}$ order.

The above examples cover all types of platonic bodies, except for the icosahedron. It was later found by A.J. Kanel-Belov that icosahedra could be put in an interlocking arrangement if a decagonal middle cross-section of the icosahedron is used as a basis[1]. The decagons can be arranged on a plane and arrows can be assigned to their faces in such a way, Fig. 9, that icosahedra reconstructed by this moving cross-section procedure produce an interlocking assembly, Fig. 10 (see the Proposition in the following section). It should be noted that not all faces of an icosahedron are in contact with its neighbors. As the dodecahedron possesses a decagonal middle cross-section as well, a similar decagon-based assembly is possible with dodecahedra, Figs. 9, 11. Thus, dodecahedra permit two interlocking arrangements normal to their symmetry axes of $3^{rd}$ and $5^{th}$ order.

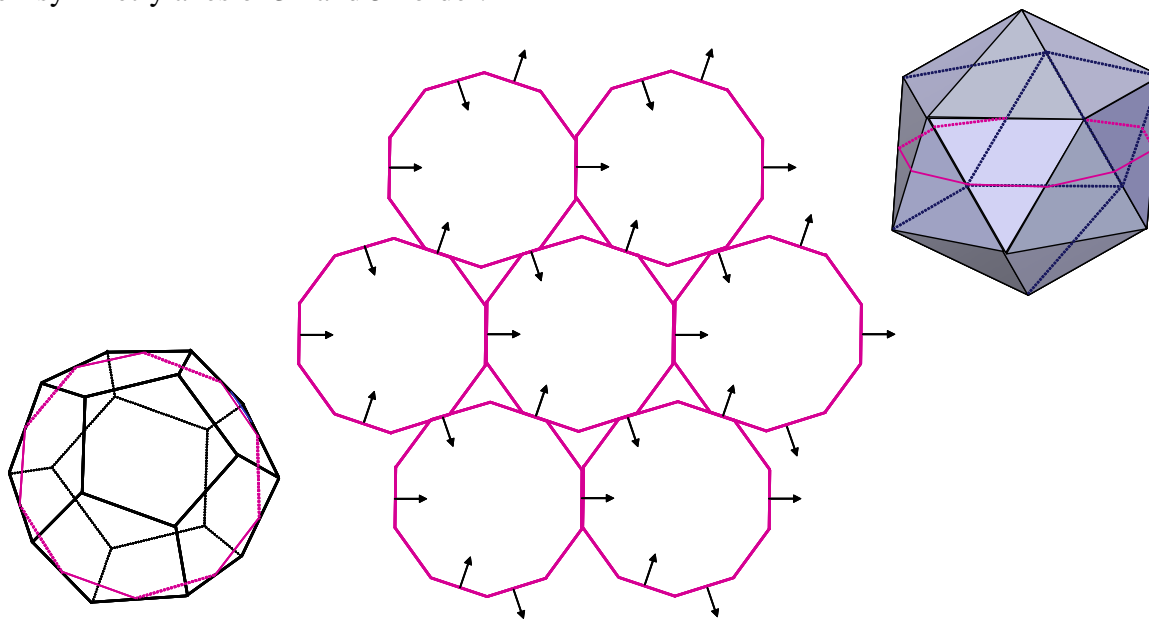

Figure 9. Decagonal tiling of plane and the associated dodecahedron and icosahedron. Only arrows that correspond to contacting faces are shown.

---

[1] It is worth mentioning that the same arrangement is also possible with truncated icosahedra (buckyballs) [7] which are the shape of $C_{60}$ molecules some fullerenes consist of.

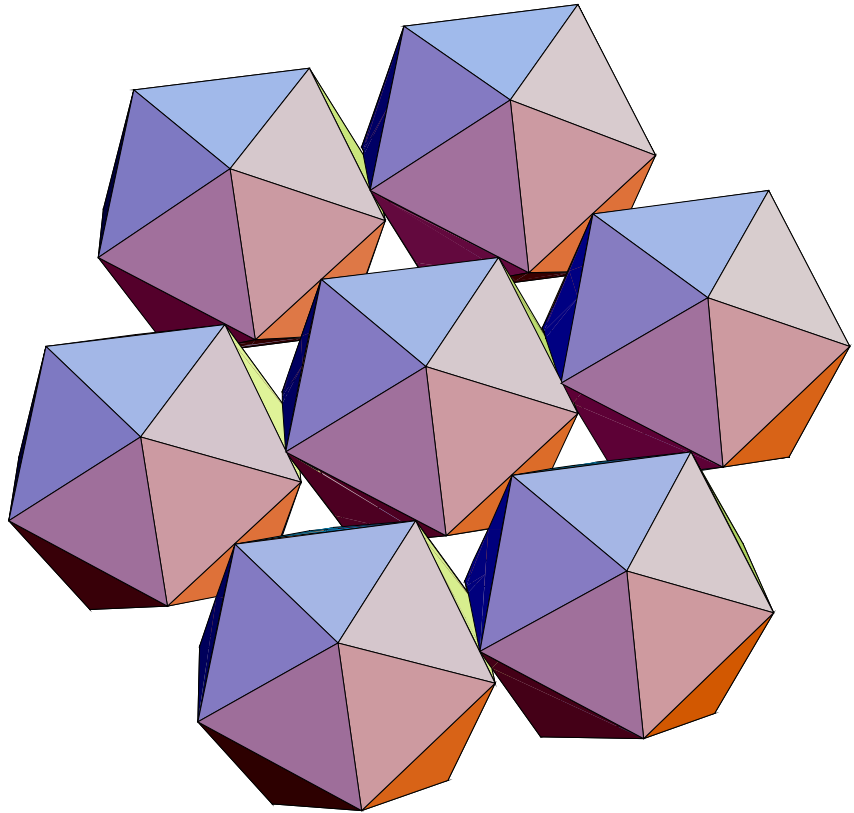

Figure 10. Fragment of interlocking set of icosahedra.

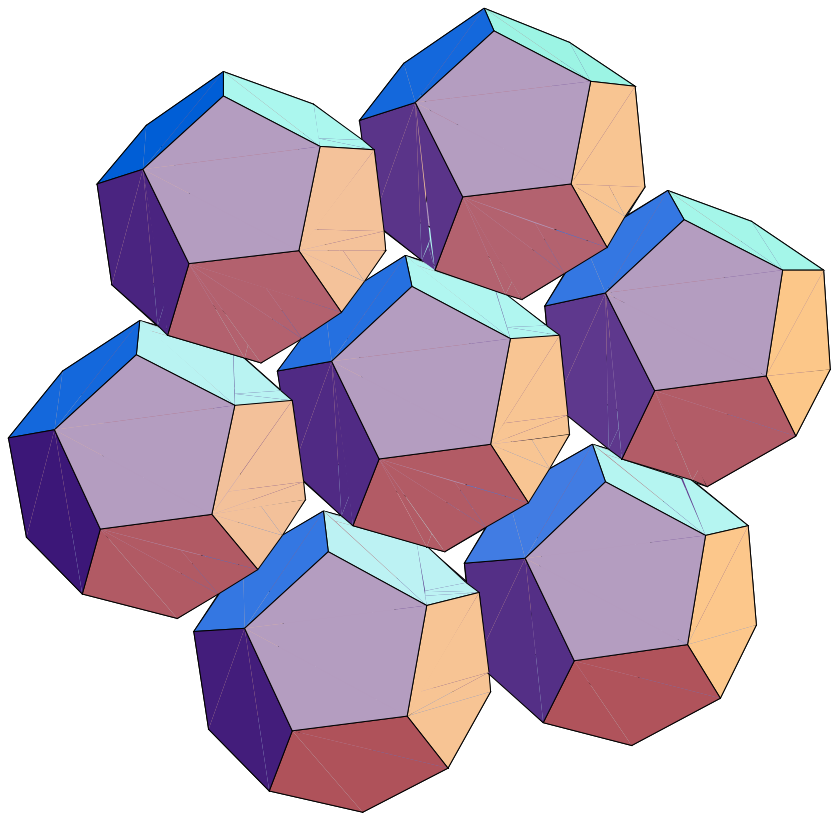

Figure 11. Fragment of interlocking set of dodecahedra normal to the symmetry axis of 5$^{th}$ order.

## 2 Geometry of interlocking: basic principles

Consider a set of three-dimensional convex solids. Let us say that two solids are *connected*, if they have contacting faces (strictly speaking: there exists a common circle on their faces). The plane containing these faces is called *a border* plane for both of these solids. Suppose that we have two solids – a solid A and a solid B respectively and their contact (border) plane. This plane defines a vector *half-space* of possible movement vectors for the solid A. Consider the set of all border planes for the given solid A.

**Proposition 1**. Let all solids in the system except A be fixed. The set of possible movement vectors of A is an intersection of all vector half-spaces of the border planes of A..

The proof is straightforward. Suppose that V is a possible movement vector. Then V belongs to all the half-spaces. On the contrary, if a vector belongs to an intersection of all vector half-spaces of the border planes, then nothing can stop the solid's movement, so V is a possible movement vector.

Suppose that we have a convex polyhedron. Consider the plane containing one face of the polyhedron. This plane divides space into two affine half-spaces. The half-space containing the polyhedron is called an *inner* half-space. The convex hull of all the affine inner half-spaces (for all faces) is obviously a polyhedron itself.

**Proposition 2.** Let all solids in the system except A be fixed. Consider contacting faces of the polyhedron A. These faces belong to border planes which define affine inner half-spaces. Let P(A) be the convex hull of these inner half-spaces. A is locked if and only if P(A) has a finite volume.

Proof. Let A be locked. Using Proposition 1 we conclude that the intersection of the border vector inner half-spaces is a zero vector. In this case border planes inscribe a finite volume: this area is a convex hull of the border affine half-spaces. Let P(A) has a finite volume. If the convex hull of affine half-spaces if finite, then the intersection of the corresponding vector spaces is a zero vector. So A cannot be moved in any direction.

Now we are ready to for discuss interlocking schemes with arrows.

**Definition 2**. The set of polyhedra is called a layer if
  1) the polyhedra are situated between two parallel planes;
  2) all of them can be produced with the finite subset using a parallel translation.

Consider a layer of solids (polyhedra) between two parallel planes. Assume that a plane P is parallel to those planes. Let the plane P be horizontal. This plane divides the polygon system into a top and a bottom part. Below we consider only the top part of the polyhedron system.

Now we mark the segments on the P plane which are situated on the border ones. In this way we obtain the faces view of an interlocking system. It is clear that every polyhedron looks as a polygon on our projection diagram. Now we assign an arrow to each of the segments using the following rules:

  1. All arrows have the same length;
  2. All arrows start in the middle of the corresponding segment and are perpendicular to it;

   3. Directions of the arrows correspond to the inclination of each border plane (the top part) relative to P.

Thus we have obtained a chessboard-like interlocking diagram with arrows, which was considered above.

Let us move the plane P upward by a small distance d. Consider the changed interlocking diagram. The segments have moved in the direction marked by an arrow by a distance d·tanα, where α is the angle between the plane P and the corresponding border plane. During this transformation the lengths of the segments may change. By varying d from zero to some value D we are executing an *evolutionary transformation* of the interlocking diagram. Let us consider a polygon subjected to this evolutionary transformation. As a result of this transformation the polygon may evolve to a segment or a point. In this case the polygon can be referred to as vanishing.

**Theorem 1. The interlocking criterion**
A layer of solids is a set with interlocking if and only if the following condition holds:
For every polyhedron the corresponding polygon vanishes during the evolutionary transformation.

Proof. It is obvious that the phase space of a polygon's evolution is exactly the region bounded by border planes of the corresponding polyhedron. So, the polygon vanishes if and only if this region is finite. Using Proposition 2 the proof is arrived at.

This theorem, along with the proof, is immediately generalized to *n* dimensions.

2.1     *Criteria of degeneration of polygons and polyhedra*

As explained above, to determine whether a solid body of a particular geometry is interlocked, it is necessary to follow the evolution, E, of its cross-section due to translations of the secant plane and to check whether the cross-section degenerates to a segment of a straight line or to a point.

**Lemma 1**: A convex polygon *M* degenerates as a result of evolution E if, and only if, one of the following occurs:
(1) The distance between a pair of parallel lines, which are continuations of some faces of *M*, vanishes.
(2) A triangle whose faces are continuations of some faces of *M* degenerates.

When interlocking in $R^4$ is considered, the following lemma is instructive.

**Lemma 2**: A convex polyhedron *M* degenerates as a result of evolution E if, and only if, one of the following occurs:
(1) The distance between a pair of parallel planes, which are continuations of some faces of *M*, vanishes.
(2) Three faces of *M* are parallel to a line and form a degenerating surface of a triangular prism.
(3) A triangular prism whose faces are continuations of some faces of *M* degenerates.

Generalisation to $\boldsymbol{R}^n$ is achieved by replacing conditions (1)-(3) in Lemma 2 with "There exist $k+1$ hyperplanes of dimension $n$-2 parallel to an $n$-1-$k$ dimensional subspace $N$ such that their projections onto the orthogonal complement of $N$ form a collapsing simplex".

2.2    *Translational interlocking and full interlocking*

In the three-dimensional structures presented in Introduction, the individual elements were shown to be interlocked, meaning that no element can be removed by translation. In principle, an element can be subjected to a combination of translations and rotations.

**Definition:** A set of solid bodies is referred to as a set with *translational interlocking* if neither of them can be removed by translation. In what follows, the term *interlocking set* will be used for brevity. If no element can be removed by any combination of translations and rotations (translational-rotational interlocking), the set is called a set with *full interlocking*. The elements of the sets will be called *translationally locked* (or simply *locked*) and *fully locked*, respectively.

Consider a vector field *v* corresponding to an infinitesimal movement of a reference element out of the assembly. An infinitesimal part of a face of a neighboring element permits the movements whose vector fields are directed inward the reference element. Therefore, *v* corresponds to an admissible movement if each infinitesimal element of neighboring faces permits it. (For a translation, *v*=const, a face blocks the movement if an infinitesimal element thereof does.) The set of vector fields corresponding to all infinitesimal movements form a 6-dimensional space in which the set of permitted infinitesimal movements form a convex cone.

We now demonstrate that the aforementioned structures are sets with full interlocking.

**Theorem 2:** Consider a layer of translationally locked polyhedra. Let $C_1$ be a minimal set of faces that ensure translational interlocking with regard to upward displacements, while $C_2$ is a minimal set of faces that ensures translational interlocking with regard to downward ones. If the polyhedra contains two spheres: $S_1$ touching faces $C_1$ and $S_2$ touching faces $C_2$, then the layer forms a set with full interlocking.

**Remark**: Such spheres do exist if (i) the faces of a polyhedron in contact with its neighbors are equally inclined to the middle plane and (ii) the circle inscribed in the middle section of the polyhedron touches its faces at points of contact with the neighboring polyhedra. This is obviously the case for the square-based and hexagonal-based arrangements shown above.

*Proof of Theorem:*

Note that an element is locked in a structure if a part of it is locked. Thus, we can use the following idea: replace each polyhedron by a sphere. We can define the sphere within the polyhedron itself so that if the sphere is locked then the whole polyhedron is locked. It is obvious that if a sphere is translationally locked then it is fully locked. So we need only to choose the sphere and check that it is translationally locked.

We will choose two different spheres: one for up-locking and another for down-locking.

Consider spheres $S_1$ and $S_2$. The centre of $S_1$ is located below the middle section, while that of $S_2$ is located above it. No movement of *P* can displace the centre of $S_1$ downward or that of $S_2$ upward, as otherwise, the spacing between them would change (which is prohibited, since rigid body motion preserves distances). The movements of the sphere centers consistent with rigid body motion are precisely those movements that are blocked by the neighboring elements. Therefore, $P$ is fully locked, and so is the element containing $P$, which completes the proof.

The above theorems establish full interlocking of assemblies shown in Figs. 2-5. Now we consider the remaining assemblies.

**Proposition:** *The layers corresponding to Fig. 7 are sets with full interlocking.* Consider Fig. 10 where the middle sections of the arrangements of icosahedra and dodecahedra normal to their symmetry axes of $5^{th}$ order are shown. Suppose that the section plane moves with a constant velocity in the direction normal to the plane of drawing such that the faces of the decagon in the section move in the directions indicated by the arrows. Due to the symmetry, the faces of the decagon will move in these directions with a constant velocity, $v$. Consider the common point of intersection of the continuations of faces $e$ and $g$ with the broken line perpendicular to face $c$. This point will move in the direction of arrow $c$ with velocity $v/\sin\varphi_e$, where $\varphi_e$ is the angle between the broken line and faces $e$ and $g$. This velocity is obviously greater than $v$, so that faces $e$ and $g$ will eventually catch up with face $c$ and the polygon will degenerate to a point. Therefore, according to Theorem 1, these layers are interlocking sets.

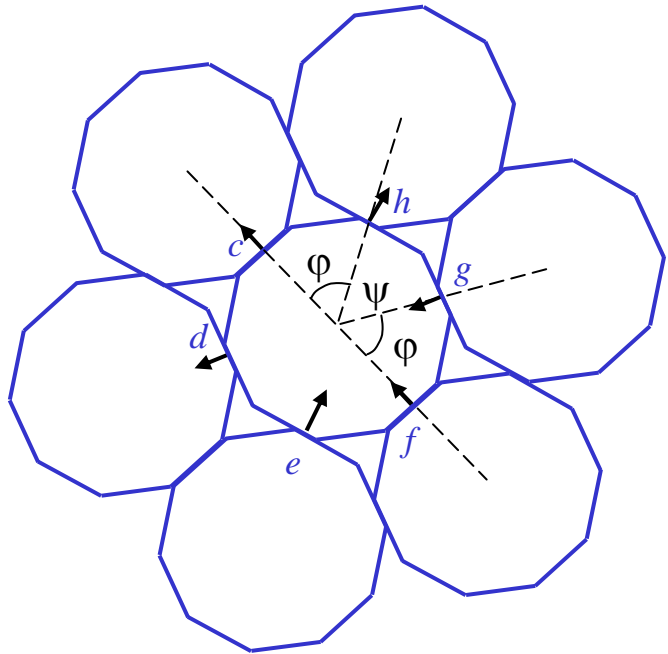

Figure 12. The contacts between the reference polyhedron and its neighbors in a decagon-based arrangement.

To prove full interlocking we note that for the assembly of icosahedra the rotation of the reference element about any axis normal to line $cf$ (Fig. 10) is blocked. Indeed, the plane contact areas at points $c$ and $f$ prevent the rotation, since each of them contains the base of the normal drawn from the centre of the reference icosahedron. Suppose there exists an axis the rotation about which is not blocked. An infinitesimal rotation about this axis is equivalent to an infinitesimal tangent vector field $v_1$. Since the normal plane passing through $cf$ is a plane of mirror symmetry, there exists a symmetrically situated axis with admissible vector field $v_2$, such that $v_1+v_2$ is a vector field parallel to the symmetry plane. This vector field corresponds to rotation about an axis normal to line *cf*, which has been shown to be blocked. This completes the proof.

The question whether in a general case a single layer with interlocking provides full interlocking remains open.

## 3   Four-dimensional cubes

In this section we generalize the recipe for generating interlocking structures of cubes to $\boldsymbol{R}^4$. The packing of three-dimensional cubes in a structure with interlocking (Fig. 4) was obtained by making use of the fact that a section of a cube by a plane normal to its principal diagonal and passing through its center is a hexagon. The evolution of a polygonal section created by the secant plane as it translates with its normal staying aligned with the principal diagonal of the cube can be described as follows. As the plane traverses the cube starting from one end of the diagonal, it first generates as its cross-section an equilateral triangle, which inflates as the plane translates. When the secant plane touches three vertices closest to the starting one, the triangle gets truncated by another equilateral triangle, which shrinks as the secant plane continues its movement. As the secant plane passes trough the cube center, these triangles become equal and the section becomes hexagonal. Further translation of the secant plane yields the same picture in an antisymmetric way up to the degeneration of the section to a point. Thus, three faces of the cross-section can be extended to form a triangle which degenerates during the evolution. According to Lemma 1, the section itself degenerates during this evolution and hence the arrangement of cubes based on the honeycomb packing of their middle sections is a set with interlocking.

A corresponding section for a four-dimensional cube is an octahedron. In analogy to the evolution of the section of a three-dimensional cube, in 4 dimensions translation of the secant hyperplane starting from a cube vertex produces a tetrahedron, which inflates as the hyperplane moves. At some moment in the evolution of the tetrahedron its vertices get truncated by another (shrinking) tetrahedron. Eventually, this results in the intersection of two identical tetrahedra - an octahedron - which is the central section of the 4D cube. Further translation of the secant hyperplane yields the same picture in an antisymmetric way up to the degeneration of the hypersection. Thus the middle section – an octahedron – can be continued to one of the tetrahedra, which degenerates as a result of the hyperplane movement.

Using the above pattern, an interlocking layer of 4D cubes can be constructed as follows. Fill $\boldsymbol{R}^3$ with octahedra by periodic translation of one octahedron in three directions. Then attach a 4D cube to each octahedron making it to its central hypersection. These cubes will be locked as, according to Lemma 2, the hypersections parallel to the middle hypersection degenerate as the secant hyperplane moves away from the middle one.

## 4   Finite interlocking structures

So far we considered infinite interlocking layers whose construction was based on regular tilings of a plane. Finite structures can be obtained by replacing the plane with, say, a sphere. Technically, it can be done as follows. Consider a sphere and match it with a sufficiently fine mesh with square and hexagonal cells. Consider the large polyhedron formed by the nodes of the mesh. Attach arrows to its edges in alternating order, in the same way as it was done for the plane. Then using each face of the large polyhedron as a base, construct a polyhedron with faces inclined as prescribed by the arrows. The polyhedra thus constructed form an interlocking envelope. (For small enough polyhedra, mutual inclinations of the adjacent faces of the large polyhedron will not violate the interlocking property.) The interior of the interlocking envelope can then be filled with arbitrary convex bodies.

The interlocking set obtained is non-convex. It can, however, be made convex if the halves of the blocks external to the large polyhedron are removed. The remaining envelope will also form a set with interlocking.

In this structure the internal convex bodies only provide a filling and it is of no importance whether they are locked in themselves, since the envelope prevents their removal from the structure. In some cases, however, it is important that the internal solids be interlocked as well.

Suppose that some elements are deleted from the layer. It is of interest to find conditions for the following property. For any finite number of missing elements in the layer a certain number of further elements can be deleted such that the rest of the layer retains interlocking.

**Definition.** A sublayer is a layer with a missing finite number of elements.

*Definition*: An interlocking layer completely filled with convex polyhedra such that each sublayer has a sublayer with interlocking is called *differential interlocking layer*.

A question whether a *differential interlocking layer* exists remains open. The existence of such a layer would mean that it is possible to construct a convex structure of convex polyhedra with interlocking, without empty space in it, such that the removal of the external layer does not lead to disintegration of the entire structure. This possibility would have far-reaching engineering applications.

## 5  Outlook

The diversity of interlocking arrangements of regular convex polyhedra uncovered so far suggests that there might be a whole wealth of such structures beyond the examples presented. These structures, as well as conventional masonry structures, cracked solids, blocky rock masses and some asteroids and comets are all examples of fragmented bodies. Considerations presented show a way to develop a geometric theory of fragmented bodies. Such a theory could also shed light on some interesting natural phenomena, e.g., the integrity and mechanics of flexible sandstones (e.g., [10]) or ancient building techniques, like dry stone walling (e.g., [11]) and may have some exciting engineering applications, such as a new concept of protective tiles for Space Shuttle [12].

**Acknowledgements**. Support from the Australian Research Council through the Discovery Grants DP0559737 (AVD) and DP0773839 (EP) is acknowledged. The first author acknowledges financial support from Israel science foundation under grant No (1178/06) and from the University of Western Australia through a Gledden Senior Visiting Fellowship. The authors are grateful to the anonymous reviewer who made a number of useful suggestions that helped improving the manuscript. We would also like to thank Anna Sorokina and Editors for professional editing of this paper.